\documentclass[12pt,psamsfonts]{amsart}
\usepackage{amssymb}
\newtheorem{Theorem}{Theorem}[section]
\newtheorem{Lemma}[Theorem]{Lemma}
\newtheorem{Proposition}[Theorem]{Proposition}
\newtheorem{Corollary}[Theorem]{Corollary}

\newtheorem{sit}[Theorem]{}

\theoremstyle{definition}

\theoremstyle{remark}
\newtheorem{remark}[Theorem]{Remark}
\newtheorem{rem}[Theorem]{Remark}
\theoremstyle{remark}

\theoremstyle{remark}

\theoremstyle{remark}

\theoremstyle{plain}

\numberwithin{equation}{section}

\newcommand{\cal}{\mathcal}

\newcommand{\Proj}{{\mathbb P}}

\newcommand{\G}{{\mathbb G}}
\newcommand{\N}{{\mathbb N}}
\newcommand{\Z}{{\mathbb Z}}

\newcommand{\C}{{\mathbb C}}

\begin{document}

\def\contentsname{\centerline {\bf Contents}}

\title{3-fold symmetric products of curves \\
as hyperbolic hypersurfaces in $\Proj^4$}

\author{Ciro Ciliberto}
\address{Dipartimento di Matematica\\
Universit\'a di Roma Tor Vergata\\ Via della Ricerca Scientifica\\
 00133 Roma, Italia}
\email{cilibert@axp.mat.uniroma2.it}

\author{Mikhail Zaidenberg}
\address{Universit\'e de Grenoble I\\ Institut Fourier\\
UMR 5582 CNRS-UJF\\ BP 74, 38402 St. Martin d'H\`eres Cedex,
France} \email{zaidenbe@ujf-grenoble.fr}

\thanks{
{\bf Acknowledgements:} This research  has been supported by the
EC research project HPRN-CT-2000-00099, EAGER. The present
collaboration started during a conference "Complex varieties" in
Oberwolfach (Aug. 1999) and continued with a visit of the first
author at the Institut Fourier of the University of Grenoble
(Dec. 1999) and a stay of the second author at the Max Planck
Institute f\"ur Mathematik in Bonn (Aug. 2002). The authors thank
all the above institutions for their support.}

\footnote{{\it 1991 Mathematics Subject Classification}:
{32Q45, 32H25, 14J70, 14H40}.\\

{\it Key words}: {Projective curve, symmetric cube, projective
3-fold, abelian variety, Kobayashi hyperbolic, deformations.}}

\date{}

\begin{abstract} We construct new examples of
Kobayashi hyperbolic hypersurfaces in $\Proj^4$. They are generic
projections of the triple symmetric product $V = C(3)$ of a
generic genus $g \ge 6$ curve $C$, smoothly embedded in $\Proj^7.$
\end{abstract}

\maketitle

\tableofcontents

\section*{Introduction}\label {?1}

\noindent

Given a generic smooth projective curve $C$ of genus $g$ over
$\C$, we consider the threefold symmetric product $T:=C(3)$ of
$C$. It is known that $T$ is Kobayashi hyperbolic if and only if
$g\geq 5$ \cite[Cor. 3.3]{sz}. Recall that by Brody's Theorem
\cite{br}, a projective algebraic variety $X$ is Kobayashi
hyperbolic (or simply {\it hyperbolic}) if there is no
non-constant holomorphic map $f:\C \to X$.

We suppose that $T$ is embedded in $\Proj^m$ for some $m$. We
notice that, since the irregularity of $T$ is $q(T):=h^1(T,{\cal
O}_T)=g$, a well known theorem of Barth-Larsen \cite[Thm. 1]{bl}
ensures that $m\geq 6$ as far as $g\ge 1$. Actually, if one
believes to Hartshorne conjecture (see \cite {lv}), then $m\geq
7$. \par

We consider the projection $\pi: \Proj^m\dashrightarrow \Proj^4$
from a general subspace $\Pi\subseteq \Proj^m$ of codimension
$5$, and we let $T'$ be the image of $T$ under $\pi$, so that
$T'$ is a hypersurface in $\Proj^4$ birational to $T$. Our main
result in this paper is the following theorem:

\begin{Theorem} \label{mainthm} If a genus $g$ curve $C$
with $g\ge 7$ is neither hyperelliptic nor trigonal and the
jacobian $J(C)$ is simple, then a generic projection $T'$ to
$\Proj^ 4$ of the threefold symmetric product $T=C(3)$ of $C$,
arbitrarily embedded into $\Proj^m$  as a threefold of degree
$d$, is a hyperbolic hypersurface in $\Proj^ 4$ of degree $d$.
The same conclusion holds for a certain special embedding
$T=C(3)\hookrightarrow \Proj^9$ in the case where $C$ is a
general plane quintic ($g=6,\,\,d=125$).
\end{Theorem}

It is known \cite{br, za} that hyperbolicity for projective
hypersurfaces is an open property. Thus small deformations of
$T'$ are again hyperbolic. Consequently, the above theorem
enables us to conclude that a very general hypersurface of degree
$d=\deg T$ in $\Proj^4$ is {\it algebraically hyperbolic}, that is
it does not contain neither rational and elliptic curves, nor
abelian surfaces (see e.g., \cite{Cl1, Cl2, Pa1, Pa2, Vo, ZR} and
the references therein for related results and optimal degree
bounds).\par

Theorem \ref{mainthm} extends to $\Proj^4$ a similar construction
from \cite{sz} of hyperbolic surfaces in $\Proj^3$ birational to
symmetric squares of curves. Concerning other explicit
constructions of Kobayashi hyperbolic projective hypersurfaces,
see e.g., \cite{sz2, sz3} and the literature therein.

The proof is divided as follows. In \S \ref {?2} we reduce the
proof of hyperbolicity of $T'$ to that of algebraic hyperbolicity
of the double surface $S'$ of $T'$. In \S \ref {?3} we show that
the only rational or elliptic curves in $S'$ could be irreducible
components of the triple curve $\Gamma'$ of $T'$. This is based
on a technical result, of independent interest, on deformations
of hyperelliptic and bielliptic curves inside abelian varieties
which generalizes an earlier one due to Pirola \cite{pi}. In \S
\ref {?4} we exclude the possibility that $\Gamma'$ contains
rational or elliptic components, thus completing the proof that
$S'$ is algebraically hyperbolic. We believe that a study of the
triple curves here deserves to be done on its own right. Finally,
in \S \ref {?5} we discuss the problem of finding, for a given
curve, projective embeddings of its symmetric products.

The lowest degree of a hyperbolic threefold $T'$ in $\Proj^4$ that
we can find with our methods is $125$, attained by the threefold
symmetric product $T$ of a general plane quintic, naturally
embedded in $\Proj^9$. Examples of lower degree can be found in
\cite {sz2}. Our approach here is mainly geometrical and exploits
only minimum of analysis. Essentially, the only analytical fact we
need is the Bloch Conjecture (see e.g., \S 9 of [De]). Notice that
the proofs in \cite {sz2} depend heavily on the value distribution
theory.
\par

\section{First part of the proof: reduction to the hyperbolicity
of the double surface}\label {?2}

 In this section we go on keeping the conventions and the
 notation introduced
above. Namely $C$ is a  general curve of genus $g\geq 5$ i.e., $C$
is neither hyperelliptic nor trigonal and $J(C)$ is simple. We
frequently use below the following elementary observation.
\par

\begin{Lemma}\label{genus} For a curve $C$ as above,
the threefold symmetric product $T=C(3)$ does not contain curves
of geometric genus $<g$.
\end{Lemma}

\begin{proof} Since $C$ is neither hyperelliptic
nor trigonal the Abel-Jacobi map
$\alpha:T=C(3)\to J(C)$ is injective. Moreover, as $J(C)$ is
supposed to be simple the image $\alpha(E)\subseteq J(C)$ of a
curve $E\subseteq T$ cannot be of geometric genus $<g$.
\end{proof}

In particular, $C$ as above does not admit a $2:1$ or a $3:1$ map
onto a curve of smaller genus.

\begin{sit} {\rm We suppose that $T=C(3)$ is embedded in
$\Proj^m$, and we let $T'$ be the general  projection of $T$ to
$\Proj^4$. We notice that $T'$ is singular. We will not describe
in full details the singularities of $T'$ here. However, if one
takes into account the description of the singularities of the
generic projection to $\Proj^3$ of a non-degenerate, smooth
surface in $\Proj^m$, $m\geq 5$ (see \cite {mo}, p.60), one can
see that:\medskip

\noindent $(i)$ ${\rm Sing}(T')$ is an irreducible surface $S'$
whose general point is a double point of $T'$ with tangent cone
formed by two distinct hyperplanes;\medskip

\noindent $(ii)$ there is a curve $\Gamma'\subseteq S'$ (possibly
reducible or empty) such that a general point of any irreducible
component of $\Gamma'$ is a triple point of $T'$ with tangent
cone formed by three independent hyperplanes. A general point of
any component of $\Gamma'$ is a triple point for $S'$, with
tangent cone formed by three independent planes meeting along the
tangent line to $\Gamma'$ at that point;\medskip

\noindent $(iii)$ there is a curve $\Delta'\subseteq S'$
(possibly reducible or empty) such that a general point of any
irreducible component of $\Delta'$ is a double point of $T'$ with
tangent cone formed by a hyperplane counted twice. A general
point of any component of $\Delta'$ is smooth for $S'$;\medskip

\noindent $(iv)$ worse singularities for $T'$ and $S'$ occur at
isolated points.\par }\end{sit}

\begin{remark} If $H$ is a general hyperplane in $\Proj^m$
containing the centre  of projection $\Pi$ then the restriction
$\pi|H$ projects generically the hyperplane section $T\cap H$ to
$\Proj^3$. If $H':=\pi(H)\subseteq\Proj^4$ then ${\rm
Sing}(\pi(T\cap H))=S'\cap H'$ is the double curve of the surface
$\pi(T\cap H)=T'\cap H'$. Furthermore, $\Gamma'\cap H'$ consists
of the triple points of $T'\cap H'$, whereas $\Delta'\cap H'$
consists of its pinch points. This yields $(i)-(iv)$. \end{remark}

\begin{sit} {\rm We denote by $S$ the pull-back of $S'$ via $\pi$.
Again we will not describe in full details $S$ and the induced
map $\pi: S\to S'$. We notice that $\pi: S\to S'$ is finite of
degree $2$, and $S$, as well as $S'$, is singular. Indeed, if
$\Gamma$ is the pull-back of $\Gamma'$ via $\pi$, the surface $S$
is singular along $\Gamma$ having, at a general point of any
irreducible component of $\Gamma$, a  double point with tangent
cone formed by two distinct planes meeting along the tangent line
to $\Gamma$ at that point. We consider the following diagram:
$$\begin{matrix}
&\tilde S&\to&\hat
S&\to&S&\hookrightarrow&T&\hookrightarrow&\Proj^m\cr
&\downarrow&&\downarrow&&\downarrow\pi&&\downarrow\pi&&\downarrow\pi\cr
&\tilde S'&\to&\hat
S'&\to&S'&\hookrightarrow&T'&\hookrightarrow&\Proj^4\cr
\end{matrix}$$
where $\hat S\to S$ [resp. $\hat S'\to S'$] is a normalization
morphism, $\tilde S\to \hat S$ [resp. $\tilde S'\to \hat S'$] is
a minimal desingularization, and the vertical arrows $\tilde S\to
\tilde S'$ and $\hat S\to \hat S'$ are induced by $\pi$. Notice
that $\pi: T\to T'$ is also a normalization morphism.\par}
\end{sit}

The following lemma reduces the proof of Theorem \ref{mainthm} to
establishing hyperbolicity of $S'$ rather than that of $T'$.\par

\begin {Lemma}\label {1streduction}
In the above setting, $T'$ is hyperbolic if and only if $S'$ is
hyperbolic.\end{Lemma}

\begin{proof} If $T'$ is hyperbolic, of course so is $S'\subseteq T'$.
Conversely, suppose $S'$ is hyperbolic whereas $T'$ is not. Let
$f: \C\to T'$ be a non-constant holomorphic map. Then $f(\C)$ is
not contained in $S'$, and therefore $f$ can be lifted to a
non-constant holomorphic map $\tilde f: \C\to T$. But since $T$ is
hyperbolic (see Corollary 3.3 of \cite {sz}), we arrive at a
contradiction.
\end{proof}

Next we reduce the proof of hyperbolicity of $S'$ to that of {\it
algebraic hyperbolicity} of $S'$. The latter means that $S'$ does
not contain neither rational nor elliptic curves \footnote{In
principle, we have to add 'and is not dominated by an abelian
surface`, but, for our purposes, we do not really need this more
restrictive definition.} or, equivalently, that every morphism
$f: E\to S'$, where $E$ is an elliptic curve, must be constant.
This is based on the following lemma.

\begin{Lemma} \label{irregularity}
The irregularity of ${\tilde S}'$ is $q({\tilde S}')\ge g\geq
5$.\end{Lemma}

\begin{proof} Let $x$ be a general point of $S'$, and let
$\pi^{-1}(x)=\{x_1,x_2\}$ with $x_1\neq x_2$. Fix an Abel-Jacobi
map $\alpha: T\to J(C)$, and consider the point
$y_x:=\alpha(x_1)+\alpha(x_2)$. Thus we have a rational map $
S'\ni x\longmapsto y_x\in J(C)$, which determines a rational map,
whence a morphism, $ \sigma: \tilde S'\to J(C)$. Since we are
assuming that $J(C)$ is simple then either $q(\tilde S')\geq g$
or $q(\tilde S')\ge g$ and $\sigma$ is constant. We denote this
constant, that depends on the choice of the centre $\Pi$ of
projection, by $c_\Pi$. Thus we have a rational map
$\G(n-5,n)\ni\Pi\longmapsto c_\Pi\in J(C)$ which in turn is
constant i.e., $c:=c_\Pi$ does not depend on the center $\Pi$ of
projection,  whence we may assume $c=0$. Let now $x_1,x_2$ be two
general points of $T$. By choosing $\Pi$ to be a general
$\Proj^{m-5}$ meeting the line $\langle x_1,x_2\rangle $, we
conclude that $\alpha(x_1)+\alpha(x_2)=0$, which is clearly
impossible. Therefore $q(\tilde S')\ge g$. \end{proof}

Now we perform our first reduction step:

\begin{Proposition} \label{2ndreduction} In the above setting,
$T'$ is hyperbolic if and only if $S'$ is algebraically
hyperbolic.\end{Proposition}

\begin {proof} By Lemma \ref {1streduction},
$T'$ is hyperbolic if and only if so is $S'$. If $S'$ is
hyperbolic then it is algebraically hyperbolic. Conversely,
suppose $S'$ is algebraically hyperbolic and assume it is not
hyperbolic. Then there is a non-constant holomorphic map $f:
\C\to S'$. If $f$ were {\it algebraically degenerate} i.e., if
the Zariski closure $Z$ of $f(\C)$ were not the whole $S'$, then
we would get a contradiction. Indeed in that case $Z$ would be a
curve of genus $0$ or $1$, contrary to the assumption that $S'$
is algebraically hyperbolic. Thus $f$ is not algebraically
degenerate, and therefore we can lift $f$ to a non-constant map
$\tilde f: \C\to \tilde S'$. Since by Lemma \ref {irregularity}
$q(\tilde S')\geq 5$, by Bloch's Theorem $\tilde f$ is
algebraically degenerate (see e.g., \S 9 of \cite {de}). This
implies that $f$ itself is algebraically degenerate, a
contradiction.\end{proof}

\section{Second part of the proof: deformations of
hyperelliptic and bielliptic curves in abelian varieties}\label
{?3}

We go on with the same hypotheses and notation as above. In this
section we make another step towards the proof of \ref{mainthm} by
proving

\begin{Proposition}\label{rational} For any rational or
elliptic curve $E$ and for any non-constant morphism $f: E\to
S'$, the image $f(E)$ must be contained in the triple curve
$\Gamma'$ of $S'$.\end{Proposition}

Our proof relies on a proposition on deformations of hyperelliptic
and bielliptic curves on an abelian variety  partially due to
Pirola \cite {pi}. We will actually prove an extension of Pirola's
result. As a matter of fact, it would not be difficult to prove
an even more general version concerning deformations of a curve
inside an abelian variety $A$ which is a double cover of another
curve \footnote{That is, which possesses a non-trivial holomorphic
involution.}, but we will not dwell on this here. A first version
of this result is contained in unpublished notes in collaboration
between G. van der Geer and the first author.

\begin{sit}\label{konts}{\rm
In order to state our result, we need to introduce some notation.
We let $Y$ be a projective variety, $X$ be a smooth, irreducible,
projective curve of genus $\gamma\ge 2$ and $f: X\to Y$ be a
morphism birational onto its image. We denote by $\beta$ the
homology class of $f(X)$, and we consider the {\it Kontsevich
space} ${\cal M}_\gamma(Y,\beta)$ (see \S 0.4 of \cite {fp}), so
that the isomorphism class of $f: X\to Y$ corresponds to a point
of ${\cal M}_\gamma(Y,\beta)$. One has an obvious forgetful map
$\phi: {\cal M}_\gamma(Y,\beta)\to {\cal M}_\gamma$ to the moduli
space ${\cal M}_\gamma$ of curves of genus $\gamma$. Suppose that
the isomorphism class of $X$ in ${\cal M}_\gamma$ lies in some
irreducible subvariety ${\cal H}$ of ${\cal M}_\gamma$. We denote
by ${\cal H}(Y,\beta)$ the scheme $\phi^*({\cal H})$; again the
class of $f: X\to Y$ corresponds to a point in ${\cal
H}(Y,\beta)$. We let $\delta_f({\cal H}(Y,\beta))$ be the
dimension of ${\cal H}(Y,\beta)$ at the point corresponding to
the class of $f: X\to Y$, that is the dimension of the
deformation space of $f$ within ${\cal H}(Y,\beta)$.
\par

We let ${\cal H}_\gamma$ be the hyperelliptic locus and ${\cal
E}_\gamma$ be the bielliptic locus in ${\cal M}_\gamma$. For
curves in abelian varieties, we have the following rigidity
result.}\end{sit}
\par

\begin{Proposition}\label {defohyperelliptic} Let $X$ be a smooth,
irreducible, projective curve of genus $\gamma$, let $A$ be an
abelian variety, and let $f: X\to A$ be a morphism which is
birational onto its image, whose homology class in $A$ is
$\beta$. One has:\par

\noindent $(i)$ if $X$ is hyperelliptic then $\delta_f({\cal
H}_\gamma(A,\beta))=\dim A$;\par

\noindent $(ii)$ if $A$ contains no elliptic curve and $X$ is
bielliptic then
$\delta_f({\cal E}_\gamma(A,\beta))\leq \dim A+1$
\footnote{This bound is sharp; see Example 1 in \cite[Sect.
1]{pi}.}.\end{Proposition}

\begin{proof} We notice that $A$ acts in a natural
way on ${\cal H}_\gamma(A,\beta)$. Indeed $a\in A$ acts on the
class of $f: X\to A$ by sending it to the class of $\tau_a\circ
f:  X\to A$, where $\tau_a: A\to A$ is the translation by $a$.
Since this action is faithful, we see that $\delta_f({\cal
H}_\gamma(A,\beta))\geq\dim A$. In order to prove $(i)$ one has
to prove that, up to the above action of $A$, there is no
non-trivial deformation of $f: X\to A$ to a family of maps of
hyperelliptic curves to $A$.
\par

Let ${\cal X}\to B$ be a deformation of $X$ over a disc $B$, such
that its general fiber is still hyperelliptic, and let $F: {\cal
X}\to A$ be a deformation of $f$. We may assume that ${\cal X}\to
B$ has a section which is a Weierstrass point of the $g^1_2$ on
each curve of the family. Using this we can change $F$ by
translations in $A$ in such a way that the image of all these
Weierstrass points is $0$ in $A$. \par

For a general point $b\in B$, we let $X_b$ be the corresponding
fiber of ${\cal X}\to B$, and we let $x_1+x_2$ be a divisor of the
$g^1_2$ on $X_b$. As $g^1_2$ is parametrized by $\Proj^1$ and $A$
does not contain rational curves, the point $F(x_1)+F(x_2)\in A$
does not depend on $x_1+x_2\in g^1_2$, whence it must be $0$.
Consequently, the Weierstrass points of the $g^1_2$ are sent by
$F$ to points of order $2$ of $A$.\par

Let $N_f$ be the normal sheaf to the map $f$, defined as the
cokernel of the differential $df: \Theta_X\to f^*\Theta_A$.
Recall (see \cite {ac}, p. 344 or \cite {ac2}, \S\S 5,6) that
there is an exact sequence
$$0\to {\cal T}\to N_f\to N'_f\to 0\,,$$
\noindent where ${\cal T}$ is a torsion sheaf supported by the
ramification divisor $R_f$ of $f$ (actually, ${\cal T}\cong {\cal
O}_{R_f}$), and $N'_f$ is locally free of rank $\dim A-1$. Let
$s\in H^0(X,N_f)$ be the section corresponding to the deformation
$F$ of $f$ \cite{ho}. This gives us a section $s'\in H^0(X,N'_f)$
which vanishes at the $2\gamma+2$ Weierstrass points of the
$g^1_2$ on $X$. Indeed, as we have seen above, their images are
fixed under the deformation $F$ of $f$. Since
$$c_1(N'_f)\leq c_1(N_f)=c_1(\omega_X)=2\gamma-2\,$$ we see that
$s'$ is identically zero. Hence $s$ vanishes at the general point
of $X$, and so it must be identically zero too. This proves
$(i)$.\par

The proof of $(ii)$ is quite similar. Indeed, we have to show
that, up to the action of $A$, there is at most a one-parameter
family of deformations of $f: X\to A$ to a family of maps of
bielliptic curves to $A$. \par

Let ${\cal X}\to B$ be a deformation of $X$ over a disc $B$ such
that its general fiber is
 bielliptic, and let $F: {\cal X}\to A$ be a non-trivial
 deformation of $f$.
 Arguing as in case
$(i)$, we may assume that ${\cal X}\to B$ has a section which is
a branch point
 of the bielliptic involution on each curve of the family.
 Therefore we can change $F$
by composing it with appropriate translations in $A$ in such a way
that the image of all these branch points is $0$ in $A$. \par

For a general point $b\in B$, we let $X_b$ be the corresponding
fiber of ${\cal X}\to B$, and we let $x_1+x_2$ be a divisor of the
bielliptic involution on $X_b$. Since $A$ does not contain
elliptic curves, the point $F(x_1)+F(x_2)\in A$ does not depend
on $x_1+x_2$ in the bielliptic involution, whence (by virtue of
the above convention) it must be $0$. Henceforth the $2\gamma-2$
branch points of the bielliptic involution are sent by $F$ to
points of order $2$ of $A$.\par

Let $s\in H^0(X,N_f)$ be a non-zero section corresponding to the
deformation $F$ of $f$ \cite{ho}.  Then the associate section
$s'\in H^0(X,N'_f)$ vanishes at $2\gamma-2$ distinct points on
$X$, which, in principle, gives no contradiction. We claim however
that any two such non-zero sections $s_1$ and $s_2$ define the
same line subbundle of $N_f$. Indeed, $N_f$ is a quotient of the
trivial bundle $f^*\Theta_A$, whence is spanned. If $\dim A=d$
then we can find sections $s_3,...,s_{d-1}\in H^0(X,N_f)$ such
that, at a general point $x\in X$, the vectors
$s_3(x),...,s_{d-1}(x)\in N_{f,x}\cong\C^{d-1}$ are linearly
independent and, moreover, the linear subspaces ${\rm
span}(s_1(x),s_{2}(x))$ and ${\rm span}(s_3(x),...,s_{d-1}(x))$
of $N_{f,x}$ are transversal.

 From the exact sequence
$$0\to
\Theta_X\stackrel{df}{\longrightarrow} f^*\Theta_A\to N_f\to 0$$
one obtains (see e.g., \cite{ha}, Ch. II, Exercises 5.16(d) and
6.11-6.12):
$$\bigwedge^{d-1}N_f=\det N_f\cong (\det \Theta_X)^{-1}=\omega_X\,.$$
On the other hand, the holomorphic section $\sigma:= s_1\wedge
s_2\wedge s_3\wedge...\wedge s_{d-1}$ of $\bigwedge^{d-1}N_f\cong
\omega_X$ vanishes to order $2$ at $2\gamma-2$ points, whence it
must be identically zero. This implies that at a general point
$x\in X$, the vectors $s_1(x), s_2(x)$ are linearly dependent,
proving our claim.\par

Now we can conclude that there is no pair $s_1, s_2$ of linearly
independent sections as above, thus proving the assertion. Indeed
for a general point $x\in X$ there exists a non-zero linear
combination $s''$ of $s_1, s_2$ which vanishes at $x$. As $s''$
also vanishes at the $2\gamma-2$ points where $s_1$ and $s_2$
vanish, $s''$ must be identically zero and so, $s_1, s_2$ must be
linearly dependent, as claimed.\end{proof}

As a consequence, we have the following

\begin{Corollary}\label{defohyperellipticcor} We let $X$ be a smooth,
irreducible, projective curve of genus $\gamma$, $A$ be an
abelian variety and $Y$ be a closed, irreducible subvariety of
$A$. If $f: X\to Y$ is a morphism birational onto its image,
whose homology class is $\beta$, then the following hold. \par

\noindent $(i)$ If $X$ is hyperelliptic then $\delta_f({\cal
H}_\gamma(Y,\beta))\leq\dim Y$.\par

\noindent $(ii)$ If $A$ contains no elliptic curve
and $X$ is
bielliptic then $\delta_f({\cal E}_\gamma(Y,\beta))\leq \dim
Y+1$.\end{Corollary}

\begin{proof}  Let us prove $(i)$;
the proof of $(ii)$ is similar and we leave it to the reader. By
Proposition \ref {defohyperelliptic}, for any point $a\in A$ the
maximal deformation family of $f$ such that the image of $X$
under the corresponding maps contains $a$, consists of just a
$1$-dimensional family of translations of $f: X\to A$. Therefore
for any point $a\in Y$, there is at most a $1$-dimensional family
of deformations of $f: X\to Y$ in ${\cal H}_\gamma(Y,\beta)$ such
that the image of $X$ under the corresponding maps passes through
$a$. This immediately implies the assertion.
\end{proof}

We are now in a position to prove \ref{rational}.\medskip

\noindent {\it Proof of Proposition \ref {rational}.} First we
exclude the existence of a morphism $f: \Proj^1\to S'$ birational
onto its image $Z'$, where $Z'\not\subseteq\Gamma'$. We suppose
that such a morphism does exist for a generic projection
$\pi:\Proj^m \dashrightarrow \Proj^4$, and we consider the curve
$Z:=\pi^{-1}(Z')$ and its normalization $f:X\to Z\subseteq T$. As
$Z'\not\subseteq \Gamma'$, the smooth curve $X$ admits a
$2$-to-$1$ morphism $\pi\circ f$ to $Z'$, whence it has an induced
$2$-to-$1$ morphism to $\Proj^1$. Since by our assumption $C$ is
neither hyperelliptic nor trigonal, the 3-fold $T=C(3)$ contains
no rational curve. Henceforth the curve $X$ is irreducible and
hyperelliptic. We show below that $\delta_f({\cal
H}_\gamma(T,\beta))\geq 5$, which contradicts part $(i)$ of
Corollary \ref {defohyperellipticcor}. Indeed, since $T=C(3)$
contains no rational curve, any Abel-Jacobi map $\alpha :
T=C(3)\to J(C)$ is injective, and so in Corollary \ref
{defohyperellipticcor} we can take $A=J(C)$ and $Y=\alpha(T)$.
\par

We may assume that $T\subseteq \Proj^7$. We chose a general plane
$\Pi\subseteq \Proj^7$ for the centre of projection to $\Proj^4$.
We consider a hyperelliptic curve $X$ and a map $f: X\to T$ as
above. For a general divisor $D=x_1+x_2$ in the $g^1_2$ on $X$,
we consider the line $\ell_D=\langle f(x_1),f(x_2)\rangle
\subseteq \Proj^7$. The Zariski closure of the union of all lines
$\ell_D$ with $D\in g^1_2\simeq \Proj^1$ is a surface scroll
$\Sigma$, and the plane $\Pi$ intersects any fiber of the induced
ruling $\Sigma\dashrightarrow \Proj^1$. Notice that $\Sigma$
cannot be a plane. Indeed otherwise $Z=f(X)\subseteq\Sigma$ would
be a conic, which is impossible, because $T$ does not contain
rational curves. Thus the maximal dimension of a family $\Im$ of
planes intersecting any fiber of the ruling of $\Sigma$ is $10$
attained in the case where $\Sigma$ is a cone and the planes in
question are passing through its vertex. When $\Pi$ runs over
$\Im$, and only for those planes $\Pi$, $Z'$ and $f:X\to Z$ do
not vary. Due to our assumption above, this clearly yields
 $$\delta_f({\cal H}_\gamma(Y,\beta))\geq \dim \G(2,7)-10=5\,,$$
 proving $(i)$.
\par

Repeating word-by-word  the above arguments and making use of
Corollary \ref {defohyperellipticcor}$(ii)$ one can exclude the
existence, for a generic projection $\pi:
\Proj^m\dashrightarrow\Proj^4$, of a morphism $f: E\to S'$ from
an elliptic curve $E$ birational onto its image
$Z'\not\subseteq\Gamma'$. This gives $(ii)$. We leave the details
to the reader. \qed

\section{Third part of the proof: hyperbolicity
of the triple curve}\label {?4}

Again we keep the same conventions as above. That is, we still
assume $C$ to be a neither hyperelliptic nor trigonal curve of
genus $g\geq 5$ with a simple jacobian $J(C)$.
\par

In this section we conclude the proof of Theorem \ref{mainthm}.
Its first claim follows from Propositions \ref{2ndreduction} and
\ref{rational} by virtue of the following result.

\begin{Proposition}\label{triplecurve} In the above setting, if in
addition $g\geq 7$, then no irreducible component of the triple
curve $\Gamma'$ of $T'$ is rational or elliptic.\end{Proposition}

Let us introduce some notation and make several useful comments.
For the time being (until Lemma \ref{4l}$(ii)$ below) it will be
sufficient to suppose $g\ge 5$.

\begin{sit}\label{s1} {\rm  Assuming
that $T\subseteq\Proj^7$, we can factor the generic projection
$\pi:\Proj^7\dashrightarrow\Proj^4$ into a projection ${\bar
\pi}:\Proj^7\dashrightarrow\Proj^6$ from a general point, which we
fix once and forever, and a projection
$\pi_L:\Proj^6\dashrightarrow \Proj^4$ from a general line
$L\subseteq\Proj^6$, which we let vary. Thus $T':=T'_L$ and
$\Gamma':=\Gamma'_L$ depend on $L$. We observe that, when we
project to $\Proj^6$, the image $\bar  T$ of $T$ acquires at worst
finitely many double points. In particular, Lemma \ref{genus}
equally applies to $\bar  T$.}\end{sit}

\begin{sit}\label{s2}{\rm  We consider the incidence relation
$$\begin{matrix}
&{\cal J}&\subseteq &T(3)&\times & \G(1,6)&\\
&&{p_1}\swarrow&&&\searrow{p_2}\\
&&T(3)&&&\,\,\,\,\,\,\,\,\,\,\,\,\G(1,6)\\
\end{matrix}
$$
where $\cal J$ is the Zariski closure of the set of pairs
$(x_1+x_2+x_3,L)$ such that $L\cap \bar T=\emptyset$ and the
points $y_i:={\bar  \pi}(x_i)\in\bar T\,\,\, (i=1,2,3)$ are
distinct, whereas $\pi_L(y_1)=\pi_L(y_2)=\pi_L(y_3)\in T'$ is a
point of $\Gamma'_{L}$. }\end{sit}

\begin{sit}\label{s4} {\rm For a general point
$\xi:=x_1+x_2+x_3\in T(3)$, we let $\Lambda_{\xi}:=\langle
y_1,y_2,y_3\rangle $ be the trisecant plane to $\bar  T$ in
$\Proj^6$ through the points $y_1,y_2,y_3$. From the General
Position Theorem \cite[p. 109]{acgh} (see also \cite[Cor.
1.3]{cc}) it follows that $\Lambda_{\xi}\cap \bar  T=\{y_1,y_2,
y_3\}$, and so the map $\xi\longmapsto \Lambda_{\xi}$ is
generically one-to-one. \par

Furthermore,
 $\,\,\pi_L(y_1)=\pi_L(y_2)=\pi_L(y_3)$ if and only
if $L\subseteq \Lambda_{\xi}$. Thus the fiber $p_1^{-1}(\xi)$ of
the first projection $p_1:{\cal J}\to T(3)$ can be naturally
identified with the dual projective plane $\Lambda_{\xi}^*\simeq
\Proj^2$. Since $T(3)$ and the general fibers $\Lambda_{\xi}^*$ of
$p_1$ are irreducible, there is only one irreducible component
${\cal J}_0$ of ${\cal J}$ which dominates $T(3)$ via the first
projection. One has $\dim {\cal J}_0=\dim T(3)+\dim
\Lambda_{\xi}^*=11$. Moreover the following holds.}
\end{sit}

\begin{Lemma} \label{1l} The map $p_2|{\cal J_0}: {\cal J_0}\to
\G(1,6)$ is surjective. Therefore also the map $p_2|{\cal J}:
{\cal J}\to \G(1,6)$ is surjective, and its fiber over a general
point $L\in \G(1,6)$ is birational to the triple curve
$\Gamma_L'$. In particular, $\Gamma_L'\neq \emptyset$.\end{Lemma}

\begin{proof} To prove the first assertion, we must
show that a general line $L$ in $\Proj^6$ is contained in a
3-secant plane $\Lambda$ to $\bar T$. To this point, we consider a
general hyperplane $H$ in $\Proj^6$ containing $L$. Clearly, the
hyperplane section $H\cap \bar T$ is a smooth linearly
non-degenerate surface in $H\simeq \Proj^5$. By a result of
Chiantini and Coppens \cite[Sect. 2]{cc}, if $L$ were not
contained in a 3-secant plane to $H\cap \bar T$ then $H\cap \bar
T$ would be either a cone or a rational normal surface of degree
4 in $H\simeq \Proj^5$. Anyhow, it would be covered by a family of
rational curves, which is excluded by Lemma \ref{genus}. This
proves the first assertion; the second one follows easily.
\end{proof}

\begin{sit}\label{ss5} {\rm By Lemma \ref{1l}, for any irreducible
component ${\cal J}'$ of ${\cal J}$ that dominates $\G(1,6)$ via
the second projection $p_2$, the fiber $p_2^{-1}(L)\cap {\cal J}'$
over a general point $L\in \G(1,6)$ is a curve birational to the
union of some irreducible components of the triple curve
$\Gamma'_L$. Thus we have $\dim {\cal J}'=\dim \G(1,6)+\dim
\Gamma'_L=11$. Furthermore, ${\cal J}'$ being irreducible, the
monodromy of the family $p_2|{\cal J}': {\cal J}'\to \G(1,6)$
acts transitively on the set of irreducible components of the
fiber $p_2^{-1}(L)\cap {\cal J}'$. Hence all these components have
the same geometric genus, which we denote by $\delta ({\cal
J}')$.}\end{sit}

\begin{sit}\label{ss6} {\rm
To prove Proposition \ref{triplecurve} we must show that $\delta
({\cal J}')\ge 2$ for any irreducible component ${\cal J}'$ of
${\cal J}$ that dominates $\G(1,6)$ via the second projection.
Arguing by contradiction, we suppose in the sequel that $\delta
({\cal J}')\leq 1$. We begin by considering a component ${\cal
J}'$ of ${\cal J}$ different from ${\cal J}_0$, assuming it does
exist. Thus ${\cal J}'$ does not dominate $T(3)$ via the first
projection $p_1$ (see \ref{s4}). }\end{sit}

\begin{sit}\label{sss6}
{\rm Clearly, for a general point $(\xi=x_1+x_2+x_3, L)$ of such a
component ${\cal J}'$, the points $y_i:=\bar\pi (x_i)\in\bar
T\,\,\,(i=1,2,3)$ are collinear (cf. \ref{s4}), and so belong to a
3-secant line $l_\xi$ to $\bar T$ that meets $L$. Thus given a
general point $\xi\in p_1({\cal J}')\subseteq T(3)$, the fiber
$p_1^{-1}(\xi)\subseteq {\cal J}'$ over $\xi$ is contained in the
set of all pairs $(\xi, L)$ such that the line $L$ meets $l_\xi$.
We denote by $G(l_\xi)$ the set of all lines $L$ in $\G(1,6)$
with $L\cap l_\xi\neq \emptyset$. We have the following
lemma.}\end{sit}

\begin{Lemma} \label{xi}
The fiber of $p_1|{\cal J}'$ over a general point $\xi$ in
$p_1({\cal J}')$ is:
$$p_1^{-1}(\xi)\cap {\cal J}'= \{\xi\}\times G(l_\xi)\,.$$
Thus $\dim p_1({\cal J}') = 11-\dim G(l_\xi)=5.$\end{Lemma}

\begin{proof} A general point $(\xi, L)\in {\cal J}'$ does
not belong to any other irreducible component of ${\cal J}$, and
the point $\pi_{L}(y_1)=\pi_{L}(y_2)=\pi_{L}(y_3)\in \Gamma'_{L}$
is smooth (see \ref{ss5}). For any line $L'\in G(l_\xi)$ we still
have $\pi_{L'}(y_1)=\pi_{L'}(y_2)=\pi_{L'}(y_3)\in \Gamma'_{L'}$
and so, clearly, $(\xi, L')$ varies in ${\cal J}'$ when $L'$
varies in $G(l_\xi)$. Notice finally that $G(l_\xi)$ is an
irreducible variety which is fibered over $l_\xi\simeq\Proj^1$
with fibers isomorphic to $\Proj^5$. This proves our
assertions.\end{proof}

\begin{sit}\label{ss61} {\rm We let again
$\xi\in p_1({\cal J}')$ be a general point with
$\xi=x_1+x_2+x_3$. As we have seen in \ref {sss6}, the points
$y_i:=\bar\pi (x_i)\in\bar T\,\,\,(i=1,2,3)$ belong to a 3-secant
line $l_\xi$ to $\bar T$. We let ${\cal V}({\cal J}')$ be the
closure in $\G(1,6)$ of the set of all these lines $l_\xi$ as
$\xi$ varies in $p_1({\cal J}')$. As $\bar T$ does not contain
lines (see Lemma \ref{genus}) the map $p_1({\cal J}')\to {\cal
V}({\cal J}')$, $\xi\longmapsto l_\xi$, is birational and has
finite fibers. Thus by Lemma \ref{xi}, ${\cal V}({\cal J}')$ is an
irreducible variety of dimension 5.

Let us also consider the variety $V({\cal J}'):=\cup_{l\in {\cal
V}({\cal J}')}l$ in $ \Proj^6$.} \end{sit}

\begin{Lemma}\label {hypersurface} $V({\cal
J}')$ is an irreducible hypersurface in $ \Proj^6$.\end{Lemma}

\begin{proof}
There is a natural incidence relation
$$\begin{matrix}
&{\cal I}={\cal I}({\cal J}')&\subseteq & {\cal V}({\cal J}')
&\times & \Proj^6 &\\
&&{q_1}\swarrow&&&\searrow{q_2}\\
&&{\cal V}({\cal J}')&&&\,\,\,\,\,\,\,\,\,\,\,\,\Proj^6\\
\end{matrix}
$$
which projects onto ${\cal V}({\cal J}')$ with irreducible
general fibers $l_\xi$. Therefore, ${\cal I}$ and also $V({\cal
J}')=q_2({\cal I})$ are irreducible. Let us show that $V({\cal
J}')$ meets a general line $L$ in $\Proj^6$, whence $\dim V({\cal
J}')\geq 5$. Indeed, as $L\in p_2({\cal J}')=\G(1,6)$ (see
\ref{ss5}-\ref{ss6}) there exists a point $(\xi=x_1+x_2+x_3,
L)\in {\cal J}'$, and so the secant line $l_{\xi}\in {\cal
V}({\cal J}')$ meets $L$.\par

To show that $\dim V({\cal J}')\leq 5$ we exploit the theory of
foci \cite{ChCi}. We consider the second projection $q_2:{\cal
V}({\cal J}') \times  \Proj^6\to \Proj^6$ and the induced
homomorphism $dq_2$ of the tangent bundles. The kernel
$T(q_2):={\rm ker} (dq_2)$ is a rank $5$ subbundle of $T({\cal
V}({\cal J}') \times \Proj^6)$ whose restriction to every fiber
$\{\tau\}\times \Proj^6$ of $q_1$ is the normal bundle of this
fiber and is clearly trivial. We pick a general point
$\tau:=l_\xi\in {\cal V}({\cal J}')$, and we let ${\cal
I}_\tau:=q_1^{-1}(\tau)\cap {\cal I}\simeq l_\xi\simeq\Proj^1$ be
the fiber  of the first projection $q_1|{\cal I}$ over $\tau$.
Then the restriction of the vector bundle $T(q_2)$ to ${\cal
I}_\tau$ identifies with the trivial bundle $({\cal
O}_{\Proj^1})^5$ (see \cite{ChCi}, (1.3)).\par

The image under $q_2$ of the fiber ${\cal I}_\tau$ is the
3-secant line $l_\xi$ of $\bar T$, with the normal bundle
$N_{l_\xi/\Proj^6}= N_{\Proj^1/\Proj^6}\simeq ({\cal
O}_{\Proj^1}(1))^5$. The homomorphism of normal bundles induced
by the projection $q_2: ({\cal J}' \times \Proj^6, {\cal
I}_\tau)\to (\Proj^6, l_\xi)$ restricts to $T(q_2)|{\cal I}_\tau$
giving the so called {\it characteristic map} $$\lambda:({\cal
O}_{\Proj^1})^5\to ({\cal O}_{\Proj^1}(1))^5\, $$ of the family of
lines ${\cal V}({\cal J}')$ at $\tau=l_\xi$. Fixing coordinates
we can write $\lambda$ via a $5\times 5$ matrix $\Phi_\lambda$ of
linear binary forms on $l_\xi\simeq\Proj^1$, called the {\it
focal matrix} of ${\cal V}({\cal J}')$ at $\tau$. We let
$F_\lambda:=\det \Phi_\lambda$; this is a binary form of degree
$5$.\par

We assume that $V({\cal J}')=\Proj^6$. Then $q_2:{\cal I}\to
\Proj^6$ is surjective and generically finite. Thus the focal
matrix $\Phi_\lambda$ is non-degenerate at a general point of
$l_\xi\simeq\Proj^1$, whence $F_\lambda$ is not identically
zero.\par

The equation $F_\lambda=0$ defines the so called {\it focal
points} of the family ${\cal V}({\cal J}')$ on $l_\xi$
(\cite{ChCi}, (1.5)). This family can also have {\it fundamental
points} and {\it cuspidal points}. Cuspidal points sitting on the
image of a general fiber ${\cal I}_t$ correspond to its singular
points (\cite{ChCi}, (1.6)). Since in our setting ${\cal
I}_t\simeq l_\xi$ is smooth, there is no such point on $l_\xi$.
Hence by Proposition 1.7 in \cite{ChCi}, the focal points of the
family ${\cal V}({\cal J}')$ on $l_\xi$ coincide with its
fundamental points, that is with the fixed points of 1-parameter
families of deformations of the line $l_\xi$ within our family
${\cal V}({\cal J}')$. These are exactly the points where the
rank of the focal matrix $\Phi_\lambda$ drops.\par

If $\xi=x_1+x_2+x_3$ then the subfamily of ${\cal V}({\cal J}')$
of lines through the point $y_i:=\bar\pi (x_i)\in l_\xi$
($i=1,2,3$) is at least $2$-dimensional. Indeed, letting
$\nu:T^3\to T(3)$ be the natural map, we consider the preimage
$P_1:=\nu^{-1}(p_1({\cal J}'))\subseteq T^3$. By Lemma \ref{xi}
it  has dimension $5$, whence the fiber of the first projection
${\rm pr}_1:P_1\to T$ over the point $x_i$ is at least
$2$-dimensional.\par

It follows that $y_i\in l_\xi$ ($i=1,2,3$) is a fundamental point
for the family ${\cal V}({\cal J}')$, where the rank of
$\Phi_\lambda$ drops at least by $2$ (cf. the proof of
Proposition 1.7 in \cite {ChCi}). We observe that the matrix
$\Phi_\lambda$ is left-right equivalent to a non-degenerate
$5\times 5$ diagonal matrix of binary forms
diag$\left(d_1(u:v),\ldots,d_5(u:v)\right)$, where $d_i|d_{i+1}$
and $d_{i+1}/d_{i}$ are the {\it invariant polynomials} of
$\Phi_\lambda$ (see e.g., \cite{Ga}, Thm. 3 in \S VI.3 or
\cite{Co}, Thm. 3 in Appendix to Ch. 6). If $y$ is a point  of
$l_\xi\simeq \Proj^1$ where the rank of $\Phi_\lambda$ drops at
least by $2$ then at least $2$ of the $d_i$ vanish at $y$, whence
$F_\lambda=\prod_{i=1}^5 d_i$ has a multiple root at $y$. In
particular, this is so for $y=y_i$ ($i=1,2,3$). Hence, as $\deg
F_\lambda=5$, this polynomial must be identically zero, a
contradiction. This proves the lemma.
\end{proof}

\begin{sit}\label{ss62} {\rm  From the previous analysis we
deduce the following. We let $L$ be a general line in $\Proj^6$
through a general point $x\in V({\cal J}')$. By Lemma
\ref{hypersurface} the general fibers of the surjection $q_1:
{\cal I}\to V({\cal J}')$ are 1-dimensional. Hence there is a
$1$-dimensional family of lines of ${\cal V}({\cal J}')$ through
$x$. This family is parametrized by some of the irreducible
components of the fiber $p_2^{-1}(L)\cap {\cal J}'$ over $L$,
that is by some components of the triple curve $\Gamma'_L$ (see
\ref{ss5}). Thus our family of all 3-secant lines to $\bar T$
passing through $x$ describes several cones with vertex $x$, say,
$\Lambda_{x,i}$ ($1\leq i\leq \nu$), over the corresponding
irreducible components of $\Gamma'_L$. These are rational or
elliptic as we suppose $\delta({\cal J}')\leq 1$ (see
\ref{ss6}).}\end {sit}

\begin{sit}\label{ss7} {\rm We fix in the sequel an Abel-Jacobi
map $\alpha: T=C(3) \to J(C)$, and we let $\alpha_3: T(3)\to
J(C)$ be the extension of $\alpha$ to the symmetric product
$T(3)$ via $$\xi=x_1+x_2+x_3\longmapsto \alpha_3 (\xi):=
\alpha(x_1)+\alpha(x_2)+\alpha(x_3)\,.$$} \end{sit}

\begin{Lemma}\label {constant} The image $W({\cal
J}'):=\alpha_3(p_1({\cal J}'))\subseteq J(C)$ has dimension
$\mu\leq 1$.\end{Lemma}

\begin{proof} We use the following notation.
For a general point $t\in W({\cal J}')$ we let ${\cal V}_t\subset
\G(1,6)$ be the family of trisecant lines $l_\xi$ to $\bar T$
corresponding to the points $\xi$ in the fiber of $\alpha_3:
p_1({\cal J}')\to W({\cal J}')$ over $t$. Thus $\dim {\cal
V}_t=5-\mu$. For an irreducible component ${\cal V}_t'$ of ${\cal
V}_t$, we let $V_t':=\cup_{l\in {\cal V}_t'} l$, so that $V_t'$ is
swept out by the $(5-\mu)$-dimensional family of lines ${\cal
V}_t'$. We note that ${\cal V}_t'\subseteq {\cal V}({\cal J}')$
and $V_t'\subseteq V({\cal J}')$.\par

If $l\in {\cal V}_t'$ is a line through a general point $x\in
V_t'$ then $l$ is contained in one of the cones $\Lambda_{x,i}\,\,
(1\leq i\leq \nu)$ with vertex $x$ (see \ref {ss62}), which we
denote by $\Lambda_{x}(l)$. As $\Lambda_{x}(l)$ is a cone over a
rational or elliptic curve (see \ref{ss62}), it is swept out by
lines $l'$ with $l' \in {\cal V}_t'$. There is henceforth a
$1$-dimensional family of lines of ${\cal V}_t'$ through a
general point of $V_t'$, and so $\dim V_t'=\dim {\cal
V}_t'=5-\mu$. Moreover, as $V_t'$ contains the cone
$\Lambda_{x}(l)$ we have $5-\mu\geq 2$.
\par

If $5-\mu=2$ then $V_t'=\Lambda_x(l)$ is a surface in $\Proj^6$
with a $2$-dimensional family of lines ${\cal V}_t'$, whence is a
plane (see Thm. 1(1) in \cite{Ro}). This plane $\Lambda_{x}(l)$
cuts out on $\bar T$ a cubic curve, which contradicts Lemma
\ref{genus}. Therefore $\mu\le 2$. \par

If $\mu= 2$ then $V_t'$ is a threefold covered by a
$3$-dimensional family of lines ${\cal V}_t'$. Hence, by a
theorem of Severi-Segre (see Thm. 1(2) in \cite {Ro}), $V_t'$ is
either a $\Proj^3$, an irreducible quadric, or a scroll in planes
\footnote{that is, there is a unique $\Proj^2$ through a general
point of $V_t'$.}.\par

If $V_t'$ is a $\Proj^3$ described by a $3$-dimensional family of
trisecant lines to $\bar T$ then $V_t'$ cuts $\bar T$ along a
cubic surface. Thus $\bar T$ must contain a rational curve, which
contradicts Lemma \ref{genus}.\par

If $V_t'$ is a quadric then the tangent space to $V_t'$ at $x$
cuts out on $V_t'$ a quadric cone $\Lambda_{x}(l)\subseteq V_t'$.
We consider the curve $\bar \Gamma:=\Lambda_{x}(l)\cap \bar T$. As
$\bar \Gamma$ is met by the lines of $\Lambda_{x}(l)$ in three
points, $\bar \Gamma$ is either a plane cubic or a space curve of
degree $6$. Anyhow, the geometric genus of $\bar\Gamma$ is at most
$4$, which again contradicts Lemma \ref{genus}.
\par

Similarly, if $V_t$ is a scroll in planes then $\Lambda_{x}(l)$ is
a plane which cuts out on $\bar T$ a cubic curve. This yields to a
contradiction as above.\par

Thus $\mu\leq 1$, as required.\end{proof}

\begin{sit}\label{s5} {\rm We turn further to considering
the component ${\cal J}'={\cal J}_0$ of ${\cal J}$ (see \ref{s4}
and \ref{ss6}). In Lemma \ref{3l} below we establish  an analog of
Lemma \ref {constant} for ${\cal J}_0$. We begin again by some
preliminary observations.\par

 For a general point $\xi=x_1+x_2+x_3\in T(3)$, the
dual $\Lambda_{\xi}^*\simeq \Proj^2$ of the 3-secant plane
$\Lambda_{\xi}:=\langle y_1,y_2,y_3\rangle\simeq \Proj^2$ to
$\bar T$ is naturally embedded in the grassmanian $\G(1,6)$ (cf.
\ref{s4}). We let
$${\cal J}_{\xi}:=p_2^{-1}(\Lambda_{\xi}^*)\cap {\cal J}_0\qquad
\mbox{ and}\qquad \Omega_{\xi}:=p_1({\cal J}_{\xi})\subset
T(3)\,.$$
 If $(\zeta,L)\in {\cal J}_0$ is a general point with
 $\zeta=z_1+z_2+z_3\in T(3)$ then
$(\zeta,L)\in {\cal J}_{\xi}$ if and only if $L=\Lambda_{\xi}\cap
\Lambda_{\zeta}$, that is iff $\pi_L\circ \bar\pi : T\to T'=T'_L$
sends both $x_i$ and $z_i$ ($i=1,2,3$) to the triple curve
$\Gamma'_L$ (see \ref{s1}). Thus $\xi$ and $\zeta$ determine the
line $L=L_{\xi,\zeta}$ in a unique way. Therefore the morphism
$p_1:{\cal J}_{\xi}\to \Omega_{\xi}$ is birational, and so $\dim
\Omega_{\xi}=\dim {\cal J}_{\xi}=3$. Furthermore, a general point
$\zeta\in T(3)$ belongs to $\Omega_{\xi}$ if and only if the
planes $\Lambda_{\xi}$ and $\Lambda_{\zeta}$ meet along a line,
that is along $L=L_{\xi,\zeta}$.}
\end{sit}

\begin{Lemma} \label{2l} $(i)$ For any pair of general points
$\zeta,\,\zeta'\in\Omega_{\xi}$, the planes $\Lambda_{\zeta}$ and
$\Lambda_{\zeta'}$ meet at a point \footnote{which is the
intersection point of the lines
$L_{\xi,\zeta}:=\Lambda_{\zeta}\cap \Lambda_{\xi}$ and
$L_{\xi,\zeta'}:=\Lambda_{\zeta'}\cap \Lambda_{\xi}$.}.\par

$(ii)$ For any 3 general points
$\zeta_1,\,\zeta_2,\,\zeta_3\in\Omega_{\xi}$, $\xi$ is an
isolated point of the intersection $\bigcap_{i=1}^3
\Omega_{\zeta_i}\,.$\end{Lemma}

\begin{proof} $(i)$ Assuming on the contrary that the
planes $\Lambda_{\zeta}$ and $\Lambda_{\zeta'}$ meet along a line
$L_{\zeta,\zeta'}$,  we let $\Pi_{\xi,\zeta}\simeq\Proj^3$ be the
3-subspace in $\Proj^6$ spanned by the planes $\Lambda_{\xi}$ and
$\Lambda_{\zeta}$. Then, clearly, for a general point $\zeta'\in
\Omega_{\xi}$ the lines $L_{\xi,\zeta'}$ and $L_{\zeta,\zeta'}$
are distinct and are contained in $\Pi_{\xi,\zeta}$, whence also
$\Lambda_{\zeta'}\subseteq \Pi_{\xi,\zeta}$. This yields a
rational map
$$\Omega_{\xi}\ni \zeta'\longmapsto \Lambda_{\zeta'}\in
\Pi_{\xi,\zeta}^*\,$$ which is generically one-to-one (see
\ref{s4}). Thus a general $\Lambda\in \Pi_{\xi,\zeta}^*$ is a
trisecant plane of $\bar  T$. It follows that $\Pi_{\xi,\zeta}$
meets $\bar  T$ along a cubic curve. This contradicts Lemma
\ref{genus}. Now $(i)$ is proven.\par

$(ii)$ The planes $\Lambda_{\zeta_i}\,\,\,(i=1,2,3)$ meet
$\Lambda_{\xi}$ along the 3 lines $L_{\xi,\zeta_i}$ in general
position. We let
$$p_{i,j}:=L_{\xi,\zeta_i}\cap L_{\xi,\zeta_j}=\Lambda_{\zeta_i}\cap
\Lambda_{\zeta_j}\subseteq \Lambda_\xi \qquad (1\leq i < j\leq
3)\,.$$ Let us pick a point $\zeta\in\bigcap_{i=1}^3
\Omega_{\zeta_i}$ close enough to $\xi$. Then $\Lambda_{\zeta}$
intersects every plane $ \Lambda_{\zeta_i}$ along a line
$L_{\zeta,\zeta_i}\subseteq \Lambda_{\zeta}$ (see \ref{s5}), and
by (i) we have
$$p_{i,j}\in L_{\zeta,\zeta_i}\cap L_{\zeta,\zeta_j}\subseteq
\Lambda_\zeta\cap\Lambda_{\zeta_i}\cap\Lambda_{\zeta_j} \qquad
(1\leq i < j\leq 3)\,.$$ Hence $\Lambda_{\zeta}$ is the plane
through the 3 non-collinear points $p_{1,2},\,p_{1,3},\,p_{2,3}$.
This holds also for $\Lambda_\xi$, therefore
$\Lambda_{\zeta}=\Lambda_{\xi}$. Thus $\Lambda_{\zeta}\cap\bar
T=\Lambda_{\xi}\cap\bar  T$, and so
 $\zeta=\xi$ (see \ref{s4}), as required.\end{proof}

\begin{sit}\label{sf6} {\rm Recall (see \ref{ss5})
that, for a general line $L\in \G(1,6)$, all irreducible
components of the fiber $p_2^{-1}(L)\cap {\cal J}_0$ have the
same geometric genus $\delta({\cal J}_0)$. To prove Proposition
\ref{triplecurve} we must show that $\delta ({\cal J}_0)\geq 2$.
Arguing by contradiction, we assume as in \ref{ss6} that $\delta
({\cal J}_0)\leq 1$.\par

For a general point $\xi\in T(3)$, we denote by $\Omega_\xi(\xi)$
the irreducible component of $\Omega_\xi$ which contains $\xi$.
Thus $\Omega_\xi(\xi)$ is birational to the irreducible component
${\cal J}_\xi(\xi)$ of ${\cal J}_\xi$ which contains the fiber
$p_1^{-1}(\xi)\cap {\cal J}_\xi=\{\xi\}\times \Lambda_\xi^*$ (cf.
\ref{s4} and \ref{s5}). Clearly, ${\cal J}_\xi(\xi)$ is saturated
by those irreducible components of the $p_2$-fibers over the
points $L\in \Lambda_\xi^*$ which meet $p_1^{-1}(\xi)$. As $\xi$
is general, there is just one such component ${\cal J}_\xi(\xi)$
of ${\cal J}_\xi$, and it has dimension 3 (see \ref{s5}). The
following analog of Lemma \ref{constant} holds.}
\end{sit}

\begin{Lemma}\label{3l}
If $\delta ({\cal J}_0)\le 1$ then, for a general point $\xi\in
T(3)$, the Abel-Jacobi map $\alpha_3: T(3)\to J(C)$ is constant on
$\Omega_\xi(\xi)$.
\end{Lemma}

\begin{proof} Letting $\xi=x_1+x_2+x_3$ we pick a general point
$\zeta=z_1+z_2+z_3$ in $\Omega_\xi(\xi)$. We claim that
$\alpha(z_1)+\alpha(z_2)+\alpha(z_3)=
\alpha(x_1)+\alpha(x_2)+\alpha(x_3)$. Indeed, we let $L\subset
\Lambda_\xi$ be a general line and $\Gamma_{(\xi,L)}$ be the
irreducible component of the fiber $p_2^{-1}(L)\cap {\cal J}_0$
containing $(\xi,L)$. It is  birational to an irreducible
component of the triple curve $\Gamma_L'$, which has geometric
genus $\delta ({\cal J}_0)\le 1$ (see \ref{sf6}). As we suppose
that the jacobian $J(C)$ is simple, the map $\alpha_3$ must be
constant along the curve $p_1(\Gamma_{(\xi,L)})$ through $\xi$.
Henceforth, for any point $\zeta=z_1+z_2+z_3\in
p_1(\Gamma_{(\xi,L)})$, we have
$\alpha(z_1)+\alpha(z_2)+\alpha(z_3)=\alpha(x_1)+\alpha(x_2)+\alpha(x_3)$.
On the other hand, if $L\subseteq \Lambda_\xi$ is a general line
and $(\zeta, L)\in \Gamma_{(\xi,L)}$ is a general point, then
$\zeta$ is a general point in $\Omega_\xi(\xi)$, which proves our
claim.\end{proof}

We need in the sequel the following\par

\begin{Lemma}\label{4l} $(i)$
For any linear series $g_9^r$ on a curve $C$ of genus $g\ge 5$
one has $r\le 4$. \par

$(ii)$ Furthermore, if $C$ is neither hyperelliptic nor trigonal
and $g\geq 7$ then $r\leq 3$. \end{Lemma}

\begin{proof} $(i)$ By the Riemann-Roch Theorem and the Serre duality,
for a divisor $D\in g_9^r$ we have:
$$r=h^0(D)-1=9-g+h^0(K_C-D)\,.$$
Thus if $D$ is non-special then $r= 9-g\le 4$ as soon as $g\ge 5$.
If $D$ is special then by Clifford's Theorem, $r\le d/2=4.5$, so
again $r\le 4$. This proves ($i$).\par

$(ii)$ Suppose now that $C$ has a $g^4_9$ and $g\geq 6$, so that
the series is special. If the $g^4_9$ has a base point then $C$
is hyperelliptic by Clifford's theorem, which is excluded by our
assumptions. Assume the series is base point free. Then it
defines a birational morphism $C\to C'$ onto a non-degenerate
curve $C'$ of degree $9$ in $\Proj^4$. The Castelnuovo genus bound
\cite[Ch. III, p. 116] {acgh} tells us that $g\le 7$. But if $g=7$
then the residual series $|K-D|$ of the $g^4_9$ with respect to
the canonical series $K:=K_C$ is a $g^1_3$, which is also
excluded.
\end{proof}

We are ready now to prove \ref{triplecurve}.
\medskip

\noindent {\it Proof of Proposition \ref{triplecurve}.} Assume
that $\delta ({\cal J}')\le 1$ for an irreducible component
${\cal J}'$ of ${\cal J}$ which dominates $\G(1,6)$ via the
projection $p_2$. There is a natural finite surjective map
$\nu_3:T(3)\to C(9)$. By abuse of notation, we still denote by
$\alpha: C(9)\to J(C)$ an Abel-Jacobi map for $C(9)$ chosen in
such a way that for $\xi=x_1+x_2+x_3\in T(3)$ one has
$\alpha(\nu_3(\xi))= \alpha_3(\xi)$ (see \ref{ss7}).

We consider first the case ${\cal J}'={\cal J}_0$. By virtue of
Lemma \ref{3l} and of Abel's Theorem we conclude that for general
points $\xi\in T(3)$ and $\zeta\in \Omega_\xi(\xi)$, all divisors
$D_\zeta:=\nu_3(\zeta)\in C(9)$ are contained in the same complete
linear system $g^r_9$ with $r\geq 3$, which we denote by ${\cal
D}_\xi$ (so ${\cal D}_\xi\simeq \Proj^r$). And also
$D_{\zeta'}\in {\cal D}_\xi$ for a general point $\zeta'\in
\Omega_\zeta(\zeta)$.
\par

We claim that, actually, one must have $r\geq 5$, which
contradicts Lemma \ref{4l}($i$).\par

We let ${\cal D}_\xi':=\nu_3^{-1}({\cal D}_\xi)\subseteq T(3)$. As
we have seen above, for a general point $\zeta\in \Omega_\xi(\xi)$
one has $\Omega_\zeta(\zeta) \subseteq {\cal D}_\xi'$. We observe
that the map $\nu_3: T(3)\to C(9)$ is unramified outside the
diagonal. Since $\xi\in T(3)$ is a general point, there is no
ramification at $\xi$, and therefore $\xi$ is a smooth point of
${\cal D}_\xi'$. By Lemma \ref{2l}$(ii)$, for any 3 general points
$\zeta_i\in \Omega_\xi(\xi)\,\,\,(i=1,2,3)$, the 3-folds
$\Omega_{\zeta_i}(\zeta_i)\subseteq {\cal D}_\xi'$ meet at $\xi$,
and $\xi$ is an isolated point of the intersection
$\bigcap_{i=1}^3 \Omega_{\zeta_i}(\zeta_i)$. Clearly, this would
be impossible if $r=\dim {\cal D}_\xi=\dim {\cal D}_\xi'$ were at
most 4. This proves the assertion in the case ${\cal J}'={\cal
J}_0$.

We turn next to the case ${\cal J}'\neq {\cal J}_0$. By Lemma
\ref {constant}, $\alpha_3$ sends the 5-dimensional subvariety
$p_1({\cal J}')$ of $T(3)$ either to a point or to a curve. In
virtue of Lemma \ref{4l}($i$) the former is impossible, whereas
in the latter case $C$ admits a linear system $g_9^4$. Since we
are assuming that $g\geq 7$ and $C$ is neither hyperelliptic nor
trigonal, this contradicts Lemma \ref{4l}($ii$). Now the proof is
completed. \qed

\begin{sit}\label {finalremark} {\rm This concludes the proof
of the first part of Theorem \ref{mainthm}. Let us point out that
the only place in the proof where one really needs the hypothesis
$g\geq 7$ (rather than $g\geq 5$) is just at the end of the
argument, in excluding the existence of a component ${\cal
J}'\neq {\cal J}_0$ with $\delta({\cal J}')\leq 1$. \par

We use this observation in the proof of Proposition \ref
{planequintic} below, which provides the second part of Theorem
\ref{mainthm}.}\end{sit}

\begin{sit}\label{planecurve}{\rm For a smooth plane curve $C$ of
degree $d$, we let ${\cal L}={\cal O}_C(1)$. In \ref{s8} below we
associate to ${\cal L}$ a line bundle ${\cal L}(3)^s$ on $C(3)$.
According to Lemma \ref{sym}$(ii)$ and Proposition \ref
{symmetric}, ${\cal L}(3)^s$ embeds $C(3)$ into $\Proj^9$ as a
smooth threefold $T$ of degree $d^3$. This map $\phi_{{\cal
L}(3)^s}: C(3)\to T\subseteq \Proj^9$ can be described more
explicitly as follows.\par

We let $\gamma: \Proj^2\to (\Proj^2)^*$ be a natural isomorphism
and, for a point $x\in \Proj^2$, we let $\gamma_x=\gamma (x)$. We
identify $\Proj^9$ with the linear system
$\Proj(H^0(\Proj^2,{\cal O}_{\Proj^2}(3)))$ of cubic curves in
$\Proj^2$, and we let $\phi: \Proj^2(3)\to\Proj^9$,
$$
x_1+x_2+x_3\longmapsto \phi(x_1+x_2+x_3):=
\gamma_{x_1}\cdot\gamma_{x_2}\cdot\gamma_{x_3}\,.$$ Since the
natural map ${\rm Sym}^3 H^0(\Proj^2,{\cal O}_{\Proj^2}(1))\to
H^0(\Proj^2,{\cal O}_{\Proj^2}(3))$ is an isomorphism, the
restriction of $\phi$ to $C(3)\subseteq \Proj^2(3)$ coincides with
$\phi_{{\cal L}(3)^s}$, whereas its restriction to the small
diagonal $\Delta=\Proj^2\subseteq \Proj^2(3)$ yields the
triple Veronese embedding of $\Proj^2$ in $\Proj^9$.
\par

If $d\geq 6$ then all hypotheses of the first part of Theorem \ref
{mainthm} are met, and therefore the generic projection $T'$ of
$T$ in $\Proj^4$ is hyperbolic. The second part of Theorem \ref
{mainthm} deals with the case $d=5$, and is provided by the
following }\end{sit}

\begin {Proposition}\label {planequintic} If $C$ is a general
plane quintic and ${\cal L}={\cal O}_C(1)$ then the line bundle
${\cal L}(3)^s$ embeds $C(3)$ in $\Proj^9$ as a smooth threefold
$T$ of degree $125$ whose generic projection $T'$ to $\Proj^4$ is
hyperbolic.\end {Proposition}

The proof is based on the following two lemmas.

\begin{Lemma}\label{g6}
We let $C$ be a smooth non-hyperelliptic curve of genus $6$. If
there exists a component ${\cal J}'\neq {\cal J}_0$ of ${\cal J}$
with $p_2({\cal J}')=\G(1,6)$ and $\delta({\cal J}')\leq 1$ then
for any linear system $g^4_9$ on $C$, any divisor $D\in g^4_9$
and any decomposition $D=x_1+x_2+x_3$ with $x_i\in C(3)$, the
points $y_i=\bar \pi(x_i)\in \bar T\,\,(i=1,2,3)$ are collinear.
\end{Lemma}

\begin{proof} We note that the special linear
systems $g^4_9$ on $C$ form just a $1$-dimensional family
$W^4_9(C)=\{|K-q|\} \quad (q\in C)$. As $C$ is not hyperelliptic,
every such series $g^4_9$ is base point free and defines a
birational embedding of $C$ onto a non-degenerate degree $9$ curve
$C'\subseteq \Proj^4$ (see the proof of Lemma \ref{4l}).
Moreover, the monodromy of this series acts on a general divisor
$D\in g^4_9$ giving the full symmetric group (see \cite {acgh},
Lemma on p. 111). That is, the ramified covering map $\varrho:
C^9\to C(9)$ is irreducible over every $g^4_9=|K-q|\subseteq
C(9)$. The map $\varrho$ admits a natural factorization :
$$C^9\stackrel{\nu_1}{\longrightarrow}
T^3\stackrel{\nu_2}{\longrightarrow}
 T(3)\stackrel{\phi_3}{\longrightarrow} C(9)\,.$$

Let us suppose that a component ${\cal J}'\neq {\cal J}_0$ with
$p_2({\cal J}')=\G(1,6)$ and $\delta({\cal J}')\leq 1$ does exist.
We keep below the notation as in the proof of Lemma
\ref{constant}. By this Lemma, the image $\alpha_3(p_1({\cal
J}'))\subseteq J(C)$ is a curve $W({\cal J}')=W^4_9(C)\simeq C$.
For a general point $q\in C=W^4_9(C)$, we consider the family
${\cal V}_q\subset \G(1,6)$ of trisecant lines $l_\xi$ to $\bar T$
corresponding to the points $\xi$ in the fiber
$P_q:=\alpha_3^{-1}(q)\cap p_1({\cal J}')$ of $\alpha_3:
p_1({\cal J}')\to W^4_9(C)\subseteq J(C)$ over $q$. It has
dimension $4$, and every irreducible component ${\cal V}_q'$ of
${\cal V}_q$ and every component $P_q'$ of the fiber $P_q$ are
also of dimension $4$. Therefore $\nu_3(P_q')=|K-q|\subseteq
C(9)$. \par

Since the monodromy action on $|K-q|$ is irreducible, the
pull-back $\tilde P_q:= \varrho^{-1}(|K-q|)\subseteq C^9$ is also
irreducible. Hence $P_q'$ coincides with the image of $\tilde P_q$
in $T(3)$. The latter shows, by the way, that $P_q=P_q'$ is
actually irreducible. As the symmetric group $S_9$ acts on
$\tilde P_q$, it follows that for any divisor $D\in g^4_9=|K-q|$
and any decomposition $D=x_1+x_2+x_3$ with $x_i\in C(3)$, the
points $y_i=\bar \pi(x_i)\in \bar T\,\,(i=1,2,3)$ belong to a
trisecant line $l_\xi$, where $\xi=x_1+x_2+x_3\in
P_q=\nu_2\circ\nu_1(\tilde P_q)$, as stated.\end{proof}

\begin{Lemma}\label{quinticconstruction}
If $C$ is a general plane quintic then, for a general $g^4_9$ on
$C$, a general divisor $D\in g^4_9$ and a decomposition
$D=x_1+x_2+x_3$ with $x_i\in C(3)$, the points $\phi(x_i)\in
T\subseteq \Proj^9\,\, (i=1,2,3)$ are not collinear, where $\phi:
C(3)\to \Proj^9$ is as in \ref{planecurve}.
\end{Lemma}

\begin{proof} Clearly, it is sufficient
to exhibit just one such smooth plane quintic and just one divisor
$D\in g^4_9$ on $C$ as above. To construct a quintic with the
required property, we take a smooth conic $C_0$, and we fix four
points $q$ and $p_i\,\, (i=1,2,3)$ on $C_0$. By Bertini's Theorem
there exists a smooth quintic $C$ cutting on $C_0$ the divisor
$3(p_1+p_2+p_3)+q$. Conversely, the divisor cut out by $C_0$ on
$C$ is $3(p_1+p_2+p_3)+q\in |K_C|=g^5_{10}$. Then the divisor
$D:=3(p_1+p_2+p_3)$ on $C$ varies in the linear system
$|K_C-q|=g^4_9$ cut out by conics through $q$. \par

We let $x_i:=3p_i\in C(3),\,\, i=1,2,3$. The corresponding points
in $\Proj^9$ lie on the Veronese image $V_3$ of $\Proj^2=\Delta$
via cubics (see \ref{planecurve}). Since $V_3$ does not contain
lines and is the intersection of quadrics through $V_3$, it does
not admit trisecant lines neither. Hence the points $\phi(x_i)\in
V_3\cap T\subseteq \Proj^9\,\, (i=1,2,3)$ cannot be collinear.
This finishes our proof.
\end{proof}

Now we can prove \ref{planequintic}.
\medskip

\noindent {\it Proof of Proposition \ref{planequintic}.} It is
well known that a general plane quintic $C$ has gonality $4$ and
the jacobian $J(C)$ is simple (e.g., see \cite{acgh}, p. 218,
\cite{ba} and \cite{cig}, Cor. (1.2)). Thereby, similarly as in
the proof of Proposition \ref{triplecurve} it follows that
$\delta({\cal J}_0)\ge 2$. Thus it suffices to exclude the
existence of a component ${\cal J}'\neq {\cal J}_0$ of ${\cal J}$
with $p_2({\cal J}')=\G(1,6)$ and $\delta({\cal J}')\leq 1$. But
for such a component the conclusion of Lemma \ref {g6} must hold,
which contradicts Lemma \ref{quinticconstruction}. Indeed, since
the points $\phi(x_i)\in T\subseteq \Proj^9\,\, (i=1,2,3)$ are
not collinear, for a generic projection $\bar \pi :\Proj^9
\supseteq T \to \bar T\subseteq \Proj^6$ the corresponding points
$y_i=\bar \pi\circ \phi(x_i)\in \bar T\,\,(i=1,2,3)$ are not
collinear either.\qed\par

\medskip

We note that the hyperbolic hypersurfaces of degree $125$ which we
have found are of the lowest degree we can construct with these
methods (cf. Lemma \ref{deg} below). See \cite {sz2} for a better
result.

\section{On projective embeddings of symmetric products}\label {?5}

In this section we address the problem of finding, for a genus $g$
curve $C$, embeddings of the symmetric product $C(n)\,\,\,\,
(n\geq 2)\,\,$ into a projective space. \par

\begin{sit}\label{s7} {\rm Let us start by recalling
some basic facts about divisor classes on $C(n)$. Given $x\in C$,
we let $\Xi_n(x)$ be the subset of $C(n)$ consisting of all
divisors $D\in C(n)$ such that $x\le D$.  Clearly, $\Xi_n (x)$ is
a reduced divisor on $C(n)$ isomorphic to $C(n-1)$. We let $\xi$
be the class of $\Xi_n (x)$ in the Neron-Severi group $NS(C(n))$.
If $B=\sum_i a_i x_i$ is a divisor on $C$, we denote by $\Xi_n
(B)$ the divisor $\sum_ia_i\Xi_n (x_i)$ on $C(n)$. This provides
a homomorphism $\Xi_n: {\rm Pic}\, (C)\to {\rm Pic}\, (C(n))$,
whose image in $NS(C(n))$ lies in $\Z\langle \xi \rangle$.
Similarly, to a line bundle ${\cal L}$ on $C$ we associate a line
bundle $\Xi_n ({\cal L})$ on $C(n)$ with $\Xi_n ({\cal L})\equiv
(\deg {\cal L})\xi$ \footnote{As usual, $\equiv$ stands for
numerical equivalence.}, called the {\it symmetrization} of ${\cal
L}$ (cf. \cite{ma}). It is defined only up to isomorphism.
\par

 The divisor $\Delta$ in
$C(n)$ given by the diagonal is divisible by $2$, because it is
the branching divisor of the natural map $C^n\to C(n)$. We let
$\delta$ be the divisor class of $\Delta\over 2$ in
$NS(C(n))$.\par

A third basic class $\theta\in NS(C(n))$ is the class of the
pull-back of a theta divisor on $J(C)$ via an Abel-Jacobi map
$\alpha: C(n)\to J(C)$. The three classes $\xi, \delta, \theta$
are related by
$$\delta=(n+g-1)\xi-\theta$$
\noindent (see e.g.,  \cite {acgh}, Prop. (5.1), p. 358 or
\cite[L. 7] {ko}).
\par

The intersection form in the submodule $\Z\langle
\xi,\theta\rangle $ of $NS(C(n))$ is given by
\begin{equation}\label{poincare}
 \xi^i\cdot\theta^{n-i}= {{g!}\over {(g-n+i)!}}\,,\qquad i=0,...,n,
\end{equation}
\noindent as dictated by Poincar\'e's formula \cite[p. 25]{acgh}
(see \cite[L. 1]{ko}).

We notice that $NS(C(n))= \Z\langle \xi,\theta\rangle $ for every
$n\in \N$ if and only if ${\rm End}(J(C))\cong \Z$.}\end{sit}

\begin{sit}\label{s8} {\rm We consider again the natural map
$C^n\to C(n)$, and we let $p_i: C^n\to C$ be the projection to the
$i$-th factor. If ${\cal L}$ is a line bundle on $C$, we denote
by ${\cal L}^{n}$ the line bundle $\bigotimes_{i=1}^np_i^*({\cal
L})$ on $C^n$. By K\"unneth's formula one has $H^0(C^n,{\cal
L}^n)\cong H^0(C,{\cal L})^{\otimes n}$. The symmetric group
$S_n$ acts on $H^0(C^n,{\cal L}^n)$, and two of the related
irreducible representations are ${\rm Sym}^n H^0(C,{\cal L})$ and
$\bigwedge^nH^0(C,{\cal L})$. If $\sigma$ is a non-zero section in
one of these subspaces then the divisor $\sigma^*(0)$ is stable
under the natural $S_n$-action on $C^n$, and so it is the
pull-back of a divisor on $C(n)$. Actually (cf. \cite{go}) there
are two line bundles on $C(n)$, which we denote by ${\cal
L}(n)^s$, respectively, ${\cal L}(n)^a$, such that

\begin{equation}\label{sections}
H^0(C(n),{\cal L}(n)^s)\cong {\rm Sym}^n H^0(C,{\cal L}),\quad
H^0(C(n),{\cal L}(n)^a)\cong \bigwedge^n H^0(C,{\cal
L})\,.\end{equation}

We need the following facts concerning the line bundle ${\cal
L}(n)^s$.}\end{sit}

\begin{Lemma}\label{sym} One has:\par

\noindent $(i)$ ${\cal L}(n)^s=\Xi_n({\cal L})$;\par

\noindent $(ii)$ $({\cal L}(n)^s)^n=(\deg{\cal L})^n$.\end{Lemma}

\begin{proof} If $s_1,...,s_n$ are $n$ sections of $\cal L$,
we denote their symmetric product by $s_1\odot ...\odot s_n$.
Given a divisor $D=x_1+...+x_n\in C(n)$, one has
\begin{equation}\label{symprod}
s_1\odot ...\odot s_n(D)=\sum_{(i_1,...,i_n)}
s_{i_1}(x_1)...s_{i_n}(x_n)\,,\end{equation} \noindent where the
sum runs over all permutations $(i_1,...,i_n)$ of the set
$\{1,...,n\}$. Thus if $A=s^*(0)=x_1+\ldots+x_k\in {\rm
Div}_k(C)$, where $s\in H^0(C,{\cal L})$, then
$$s^{\odot
n}(y_1+\ldots +y_n)=n!s(y_1)\ldots s(y_n)=0\Leftrightarrow
y_i=x_j$$ for some $ i\in\{1,\ldots,n\}$ and
$j\in\{1,\ldots,k\}$. Hence $(s^{\odot n})^*(0)=\Xi_n(A)\in {\rm
Div}(C(n))$. This clearly implies $(i)$. Part $(ii)$ follows from
$(i)$ and the fact that $\xi^n=1$ (see (\ref{poincare})).
\end{proof}

The following result extends Lemma 3.10 of \cite {sz}.

\begin{Proposition}\label{symmetric} ${\cal L}$
is very ample on $C$ if and only if ${\cal L}(n)^s$ is very ample
on $C(n)$. \end{Proposition}

\begin{proof} Let us show that, if the line bundle ${\cal
L}(n)^s$ on $C(n)$ is very ample, then so is ${\cal L}$ as well.
This trivially holds if $n=1$. Then we proceed by induction, and
we restrict ${\cal L}(n)^s$ to $\Xi_n(x)\simeq C(n-1)$. Since,
obviously, ${\cal L}(n)^s$ restricts to ${\cal L}(n-1)^s$, we
conclude the induction.\par

To prove the converse, we assume that ${\cal L}$  is very ample on
$C$. Let us show first that ${\cal L}(n)^s$ separates points of
$C(n)$. We let $$D=p_0x_0+p_1x_1+...+p_kx_k\qquad\mbox{and}\qquad
D'=p_0'x_0+p_1'x_1+...+p_k'x_k$$ be two distinct effective
divisors of degree $n$ on $C$ regarded as points of $C(n)$, where
$x_1,\ldots,x_k$ are distinct points of $C$ and $p_i,p_i'$ are
non-negative integers with $0\le p_0<p_0'$. By the very ampleness
of ${\cal L}$, we can find two sections $s_0,\,s_1\in H^0(C,{\cal
L})$ such that $s_0(x_0)\neq 0$ and $s_1(x_0)=0$, whereas
$s_1(x_i)\neq 0$ for all $i=1,...,k$. We let $\sigma:=s_0^{\odot
p_0}\odot s_1^{\odot (n-p_0)}\in H^0(C(n),{\cal L}(n)^s)$. Then
by (\ref{symprod}) we obtain:
$$\sigma(D)=k! s_0^{p_0}(x_0)s_1^{p_1}(x_1)\cdot\ldots\cdot
s_1^{p_k}(x_k)\neq 0\,$$ whereas $\sigma(D')=0$, as $p_0'>p_0$
and $s_1(x_0)=0$. Thus $\sigma$ separates $D$ and $D'$.\par

Next we show that ${\cal L}(n)^s$ separates tangent vectors. Let
$D=x_1+...+x_n$ be a point of $C(n)$ and $X\in T_D(C(n))$ be a
non-zero tangent vector. We need to find a section $\sigma$ of
${\cal L}(n)^s$ such that $\sigma(D)=0$ whereas $X(\sigma)\not=0$.
If $D$ is formed by $n$ distinct points, the argument is similar
to those developed in \cite {sz}, {\it loc.cit}., and we will not
repeat it now. To the other extreme, suppose that $D=nx$ sits on
the small diagonal. We will prove the statement in this case only,
since the intermediate cases, where only some of the points of
$D$ come together, can be treated in a similar way.\par

Let $t$ be a local coordinate on $C$ in a small disc $W$ around
$x$. Then near $(x,...,x)$, $C^n$ is isomorphic to the polydisc
$W^n$ with coordinates, say, $(t_1,...,t_n)$. Whereas $C(n)$ about
$nx$ looks like $V:=W(n)$, which is still a polydisk with
coordinates $(z_1,...,z_n)$, where $z_i$ is the $i$-th symmetric
function on $(t_1,...,t_n)$. \par

Take now sections $s_i\in H^0(C,{\cal L})\,\,\, (i=1,...,n)$
which have the following expressions  in $W$ :
$$s_1= t+o(t),\qquad s_i=1+a_{i}t+o(t),\qquad i=2,...,n\,.$$
Since ${\cal L}$ is very ample we can find sections as above with
arbitrary values of $a_i,\,\,\, i=2,...,n$. Our objective is to
show that one can find a section of the form $\sigma=s_1\odot
...\odot s_n$ such that in $V=W(n)$ it looks like
\begin{equation} \label {esse}
\sigma=b_1z_1+...+b_nz_n+(\mbox{higher order terms in}\,\,\,
z_1,\ldots,z_n)\end{equation} \noindent with arbitrary $b_i,\,\,i
=1,...,n$. This will prove the assertion.\par

First of all, if we let $a_i=0\,\,( i=2,...,n)$ then
$$\sigma=(n-1)!z_1+(\mbox{higher order terms})\,.$$
Thus referring to (\ref {esse}) above, we can find  a section
corresponding to the $n$-tuple $(b_1,...,b_n)=(1,0,...,0)$. Let
us show by induction that for every $n$-tuple
$(0,...,0,1_j,0,...,0)$, where $j=2,...,n$, there exist
corresponding sections. By induction it suffices to prove that
for any $j=2,...,n$ we can find a section corresponding to an
$n$-tuple of the form $(b_1,...,b_j,0,...,0)$ with $b_j\not=0$.
For this it suffices to take $a_i=1,\,\,\, i=2,...,j,\,\,\,
a_i=0,\,\,\, i=j+1,...,n$, since then
$$\sigma=b_1z_1+...+b_{j-1}z_{j-1}+j!(n-j)!z_j+(\mbox{higher
order terms})$$ \noindent with some $b_1,...,b_{j-1}\in\Z$, which
we do not care to compute. This ends our proof.\end{proof}

\begin {sit}\label {antisym} {\rm In the reminder of this section
we turn to line bundles of the form ${\cal L}(n)^a$. The following
analog of Lemma \ref {sym} holds.}\end{sit}

\begin{Lemma}\label{wedge}

 $(i)$ The class of ${\cal L}(n)^a$ in $NS(C(n))$ is
$$(\deg {\cal L}-g-n+1)\xi+\theta=(\deg {\cal
L})\xi-\delta\,.$$ \noindent $(ii)$ Consequently,
$$({\cal L}(n)^a)^n=\sum_{i=0}^n
{n\choose i}{{g!}\over {(g-n+i)!}}(\deg {\cal
L}-g-n+1)^i\,\,.$$\end{Lemma}

\begin{proof} We give a proof only in the case where the line bundle
${\cal L}(n)^a$ has enough sections (cf. Remark \ref{mattuck}
below for an alternative approach). By virtue of (\ref{sections})
we may regard the sections of ${\cal L}(n)^a$ as those of
$\bigwedge^n H^0(C,{\cal L})$. Suppose that $\sigma=s_1\wedge
...\wedge s_n\in\bigwedge^n H^0(C,{\cal L})$ is such that the
divisor $\sigma^*(0)$ on $C(n)$ is reduced, where $s_1,...,s_n$
are linearly independent sections of $\cal L$. Notice that the
class of $\sigma^*(0)$ in the group NS$(C(n))$ does not depend on
the choice of $\sigma$. Given $D=x_1+...+x_n\in C(n)\backslash
\Delta$, one has
$$\sigma(D)=s_1\wedge ...\wedge s_n(D)
=\det(s_i(x_j))_{i,j=1,...,n}=0$$ if and only if there is a
non-trivial linear combination $s$ of $s_1,...,s_n$ such that
$s(x_i)=0$ for all $i=1,...,n$, that is $D\le s^*(0)$. It is not
difficult to extend this remark to the case $D\in \Delta$. That
is, also in this case $D$ belongs to $\sigma^*(0)$ if and only if
$D\le s^*(0)$ for a certain non-zero section $s\in V:={\rm
span}\langle s_1,\ldots,s_n\rangle$. We leave the details to the
reader. Thus $\sigma^*(0)=\Gamma_n(V)$, where $\Gamma_n(V)$
denotes the set of all $D\in C(n)$ subordinated to the members of
the linear system $|V|$. Now $(i)$ follows by Lemma (3.2) in
\cite[p. 342]{acgh}, whereas $(ii)$ follows from $(i)$ by the
Poincar\'e Formula (\ref{poincare}).\end{proof}

\begin{rem} \label{mattuck} Let $p:C\times C(n-1)\to C$ be the
first projection and $\alpha:C\times C(n-1)\to C(n)$ be the map
sending $(x,\sum_{i=1}^{n-1} x_i)$ to $x+\sum_{i=1}^{n-1} x_i$.
For a line bundle ${\cal L}$ over $C$, its {\it Mattuck
symmetrization} \cite[Sect. 1]{ma} is defined by the push-pull
formula :
$${\cal E}_{n, {\cal L}}:=\alpha_*p^*{\cal L}\,.$$
This is a rank $n$ vector bundle over $C(n)$. There is a canonical
isomorphism \cite[Corollary of Prop. 2]{ma} :
$$H^0(C(n), {\cal E}_{n, {\cal L}})\cong H^0(C,{\cal L})\,.$$
Indeed, for $s\in H^0(C,{\cal L})$, $p^*s$ is constant along the
fibers of $\alpha$. Also, for any open subset $U\subseteq C(n)$
one has \cite[Prop. 1]{ma} :
\begin{equation}\label{mat} H^0(U, {\cal E}_{n, {\cal L}})\cong
\bigl[H^0(\varrho^{-1}(U), \bigoplus_{i=1}^n p_i^*{\cal
L})\bigr]^{S_n}\,,\end{equation} where as above $\varrho : C^n\to
C(n)$ stands for the orbit map of the action of the symmetric
group $S_n$ on $C^n$.

Similarly, if $\Sigma_n\subseteq C\times C(n)$ is the natural
incidence relation then the {\it G\"ottsche symmetrization} of
${\cal L}$ is defined as follows \cite{go} :
$${\cal E}'_{n, {\cal L}}:=(p_2)_*(p_1)^*{\cal L}\,,$$ where
$p_1:\Sigma_n\to C$ and $p_2:\Sigma_n\to C(n)$ are the canonical
projections\footnote{Actually, \cite{go} is dealt with projective
surfaces rather than with curves.}. The map
$$\Sigma_n\ni \bigl(x,x+\sum_{i=1}^{n-1} x_i\bigr)\longmapsto
\bigl(x,\sum_{i=1}^{n-1} x_i\bigr)\in C\times C(n-1)$$ yields an
isomorphism $\Sigma_n\simeq C\times C(n-1)$ which agrees with the
above projections. This gives rise to a natural identification
$${\cal E}_{n, {\cal L}}={\cal E}'_{n, {\cal L}}\,.$$
By \cite[Thm. A.1]{go} one has: \begin{equation}\label{go1} {\cal
L}(n)^a\cong\det {\cal E}_{n, {\cal L}} \cong {\cal L}(n)^s
\otimes {\cal O}_{C(n)}(-\delta)\,,\end{equation} which by virtue
of \ref{sym}$(i)$ provides another proof of \ref {wedge}$(i)$.
Moreover (\ref{mat}) and (\ref{go1}) yield the second part of
(\ref{sections}).

In particular, if ${\cal L}=\omega_C=\Omega^1_C$ then \cite[Sect.
1, (2) ]{ma} ${\cal E}_{n, {\cal L}}\cong \Omega^1_{C(n)}$,
whence by (\ref{go1}) we have
\begin{equation}\label{canlb} \omega_{C(n)}\cong
{\omega}_C(n)^a\,.\end{equation}
\end{rem}

\begin{sit}\label{varamp}{\rm
We notice that a direct analog of Proposition \ref {symmetric}
does not hold in general for the bundles ${\cal L}(n)^a$. For
instance, if $C$ is a general plane quartic then ${\cal
L}:=\omega_C={\cal O}_C(1)$ is very ample, although ${\cal
L}(2)^a\cong \omega_{C(2)}$ is ample but not very ample (see
\cite{ko} or \cite[(4.16.b$'$)]{sz}).\par

Suppose however that ${\cal L}\to C$ is $n$-very ample i.e., any
subscheme of length $n+1$ imposes independent conditions to the
sections in $H^0(C,{\cal L})$ (see \cite{bs}). Then
$$\phi_{\cal L}: C\to \Proj( H^0(C,{\cal L})^*)$$
\noindent is an embedding, and there is no $(n+1)$-secant
$\Proj^{n-1}$ to $C$ in this embedding. The map
$$\phi_{{\cal L}(n)^a}: C(n)\to
\Proj\left(\bigwedge^n H^0(C,{\cal L})^*\right)$$
\noindent can be
described as follows. Given a divisor $D\in C(n)$, we regard it
as a subscheme of $C\subseteq \Proj( H^0(C,{\cal L})^*)$. Then $D$
spans a $\Proj^{n-1}$ which we denote by $\langle D\rangle $, and
we let
$$\phi_{{\cal L}(n)^a}(D):=\langle  D\rangle \in \G(n-1,\Proj
( H^0(C,{\cal L})^*))\subseteq \Proj\left(\bigwedge^n H^0(C,{\cal
L})^*\right).$$ As no $\langle  D\rangle\simeq \Proj^{n-1}$ is
$(n+1)$-secant to $C$, the map $\phi_{{\cal L}(n)^a}$ is
injective. Actually by \cite{cg}, $\phi_{{\cal L}(n)^a}$ is an
embedding i.e., ${\cal L}(n)^a$ is very ample, if and only if
$\cal L$ is $n$-very ample \footnote{Notice that $\phi_{{\cal
L}(n)^a}$ and the map $\varphi_{n-1}$ considered in \cite{cg} are
dual to each other.}.\par

Anyhow, the minimal degree of projective embeddings of $C(3)$
given by line bundles of the form ${\cal L}(3)^a$ is higher than
those provided by line bundles of the form ${\cal L}(3)^s$.
Indeed, we have the following lemma.}
\end{sit}

\begin{Lemma}\label{deg} \par

If a curve $C$ is neither hyperelliptic nor trigonal then for any
$3$-ample line bundle ${\cal L}$ on $C$ one has $({\cal
L}(3)^a)^3>125$.
\end{Lemma}

\begin{proof} If $\deg {\cal L}=d$ then according to
 Lemma \ref {wedge}$(ii)$ we have

\begin{equation}\label{ocen} ({\cal L}(3)^a)^3=g(g-1)(g-2) +
3g(g-1)(d-g-2)\end{equation}
$$+3g(d-g-2)^2+(d-g-2)^3.$$
\noindent In view of (\ref{ocen}) and the classification of space
curves \cite[IV.6.4.2]{ha}, it is easily seen that $({\cal
L}(3)^a)^3>125$ unless, maybe, for one of the following pairs :
\begin{equation}\label{poss}
(g,d)=(5,7),\,\,(6,5),\,\,(6,7),\,\,(6,8),\,\,
(7,8)\quad\mbox{or}\quad (8,8)\,.\end{equation} By the
 result of Catanese and G\"otsche cited in \ref{varamp}, the
bundle $({\cal L}(3)^a)^3$ is very ample and so, there is no
$4$-secant plane to the image $\varphi_{|{\cal L}|}(C)$ of $C$
under the embedding defined by the linear system $|{\cal L}|$. In
particular, we must have $\dim |{\cal L}|\ge 4$.

On the other hand, $d-g \le 2$ in all cases of (\ref{poss}).
Therefore the linear system $|{\cal L}|$ is special and different
from the canonical one. Since $C$ is neither hyperelliptic nor
trigonal, Clifford's Theorem gives the strict inequality $\dim
|{\cal L}| < d/2$, which yields $d\ge 9$. This excludes the six
remaining possibilities as in (\ref{poss}).
\end{proof}

\end{document}